\documentstyle[12pt,twoside]{article}
\newtheorem{theorem}{\bf Theorem}[section]
\newtheorem{proposition}[theorem]{\bf Proposition}
\newtheorem{lemma}[theorem]{\bf Lemma}
\newtheorem{corollary}[theorem]{\bf Corollary}
\newtheorem{example}[theorem]{\bf Example}

\newtheorem{remark}[theorem]{\bf Remark}

\input{amssym}

\date{}
\begin{document}

\title{{\Large\bf Homological and cohomological properties of Banach algebras and their second duals}}

\author{{\normalsize\sc M. J. Mehdipour and A. Rejali\footnote{Corresponding author}}}
\maketitle

{\footnotesize  {\bf Abstract.} In this paper, we investigate homological properties of Banach algebras. We show that retractions Banach algebras preserve biprojectivity, contractibility and biflatness. We also prove that contractibility of second dual of a Banach algebra implies contractibility of the Banach algebra. For a Banach algebra $A$ with $\Delta(A)\neq\emptyset$, let $\frak{F}(X, A)$ be one of the Banach algebras $C_b(X, A)$, $C_0(X, A)$, $\hbox{Lip}_\alpha(X, A)$ or $\hbox{lip}_\alpha(X, A)$. In the following, we study homological properties of Banach algebra $\frak{F}(X, A)$, especially contractibility of it.
We prove that contractibility of $\frak{F}(X, A)$ is equivalent to finiteness of $X$ and contractibility of $A$. In the case where, $A$ is commutative, we show that $\frak{F}(X, A)$ is contractible if and only if $A$ is a $C^*-$algebra and both $X$ and $\Delta(A)$ are finite. In particular, $\hbox{lip}_\alpha^0(X, A)$ is contractible if and only if $X$ is finite. We also investigate contractibility of $L^1(G, A)$ and establish $L^1(G, A)$ is contractible if and only if $G$ finite and $A$ is contractible. Finally, we show that biprojectivity of the Beurling algebra $L^1(G, \omega)$ is equivalent to compactness of $G$, however, biprojectivity of the Banach algebras $L^1(G, \omega)^{**}$ is equivalent to finiteness of $G$.  This result holds for the Banach algebra $M(G, \omega)^{**}$ instead of $L^1(G, \omega)^{**}$.}
{\footnotetext{ 2020 {\it Mathematics Subject Classification}:
 46H20, 46M20, 46M10, 46H25, 46E40, 43A20.

{\it Keywords}: Biprojectivity, contractibility, biflatness,  Lipschitz algebra, vector-valued functions, weighted group algebras.}}

\section{\normalsize\bf Introduction}

Let $A$ be a Banach algebra and $E$ be a Banach $A-$bimodule. A bounded linear operator $D: A\rightarrow E$ is called a \emph{derivation} if 
$$
D(ab)=D(a)b+aD(b)
$$
for all $a, b\in A$. Let us recall that $A$ is called \emph{contractible} if every continuous derivation $D$ from $A$ into Banach $A-$bimodule $E$ is inner, that is, there exists $z\in E$ such that $D(a)=za-az$ for all $a\in A$.
Also, $A$ is said to be \emph{weakly amenable} if every continuous derivation $D$ from $A$ into $A^*$ is inner.

For Banach spaces $E$ and $F$, a bounded linear operator $T: E\rightarrow F$ with closed range is said to be \emph{admissible} if $\hbox{ker}T$ and $T(E)$ are complemented subspaces of $E$ and $F$, respectively. Now, let $A$ be a Banach algebra and let $E$ and $F$ be Banach left $A-$modules.
We reserve the symbol $ _A{\cal B}(E, F)$ for the set of all bounded left $A-$module homomorphisms from $E$ into $F$.
A Banach left $A-$module $P$ is called \emph{projective} (respectively, \emph{flat}) if for every Banach left $A-$module $E$ and $F$, for each admissible epimorphism (respectively, monomorphism) $T\in _A{\cal B}(E, F)$ and each $S_1\in _A{\cal B}(P, F)$ (respectively, $S_2\in _A{\cal B}(E, P^*)$), there exists $R_1\in_A{\cal B}(P, E)$ (respectively, $R_2\in _A{\cal B}(F, P^*)$ ) such that $$T\circ R_1=S_1\quad (\hbox{respectively}, R_2\circ T=S_2).$$ A Banach algebra $A$ is called\emph{ biprojective} (respectively, \emph{biflat}) if $A$ is projective (respectively, flat) as a Banach $A-$bimodule.
These definitions due to Helemskii \cite{h3}. It is proved that $A$ is biprojective if and only if the corresponding diagonal operator $$\pi_A: A\hat{\otimes}A\rightarrow A$$ has a bounded bimodule right inverse. Also, $A$ is biflat if and only if $\pi_A^*$ has a bounded bimodule left inverse; for more details see \cite{d, dp, hel, hel1, h3, hel2, r1, w1}.

Johnson \cite{joh} studied the space ${\cal H}^1(L^1(G), E)$, the first Hochschild cohomology group of $L^1(G)$ with coefficients in Banach $L^1(G)-$bimodule $E$.
 He showed that ${\cal H}^1(L^1(G), E^*)$ vanishes for all Banach $L^1(G)-$bimodules $E$ if and only if $G$ is amenable. This result leads to  Banach algebras with the property ${\cal H}^1(A, E^*)$ vanishes for all Banach $A-$bimodule $E$ are called \emph{amenable}. Some authors gave equivalent definitions of amenable Banach algebras. For example, Helemskii \cite{hel, hel1, h3, h4, hel2} introduced the homological properties Banach algebras and found equivalent definitions of amenability. For instance, he showed that $A$ is an amenable Banach algebra if and only if the unitization of $A$ is biflat. For other equivalent definitions see \cite{j2}; see also \cite{h4, ram}, see Selivanov \cite{sem1} for characterization  the biprojective and biflat Banach algebra through the first Hochschild cohomology groups, \cite{sem2} for biprojective semisimple Banach algebras, for biprojectivity and biflatness of Banach algebra related to locally compact groups and semi groups see \cite{ch, gm, groh, ps, ram1, run}.

Let $X$ be a locally compact space and $A$ be a Banach algebra. Let $X^A$ be the set of all vector-valued functions from $X$ into $A$. Let us recall that $f\in X^A$ is called \emph{bounded} if
$$
\|f\|_{\infty}:=\sup_{x\in X}\|f(x)\|
$$
is finite. We denote by $C_b(X, A)$ the space of all bounded continuous functions from $X$ into $A$. We also denote by $C_0(X, A)$ the set of all continuous functions from $X$ into $A$ that vanish at infinity. If $A=\Bbb{C}$, then we  use from the symbols $C_b(X)$ and $C_0(X)$ instead of $C_b(X, A)$ and $C_0(X, A)$, respectively. 

Let $f$ be a bounded vector-valued function from metric space $(X, d)$ into Banach algebra $(A, \|.\|)$. Then $f$ is called \emph{vector-valued Lipschitz function of order} $\alpha$ if there exists $M>0$ such that for every $x, y\in X$
$$
\|f(x)-f(y)\|\leq Md^\alpha(x, y),
$$
where $\alpha$ is a positive real number. The space of all vector-valued Lipschitz functions of order $\alpha$ from $(X, d)$ into $(A, \|.\|)$ is denoted by $\hbox{Lip}_\alpha(X, A)$.   Clearly, every vector-valued Lipschitz function of order $\alpha$ is continuous. Hence $\hbox{Lip}_\alpha(X, A)$ is a subspace of $C_b(X, A)$. It is proved that $\hbox{Lip}_\alpha(X, A)$ with the pointwise multiplication and the norm
$$
\|f\|_{\infty, \alpha}:=\|f\|_{\infty}+p_\alpha(f)
$$
is a Banach algebra, where
$$
p_\alpha(f)=\sup_{x, y\in X,x\neq y}\frac{\|f(x)-f(y)\|}{d^\alpha(x, y)}.
$$
Let $\hbox{lip}_\alpha(X, A)$ be the subalgebra of $\hbox{Lip}_\alpha(X, A)$ consisting of $f$ such that
$$
\lim_{d(x, y)\rightarrow 0}\frac{\|f(x)-f(y)\|}{d(x, y)^\alpha}=0.
$$
If $X$ is a locally compact metric space, then we set $$\hbox{lip}^\circ_\alpha(X, A):=\hbox{lip}_\alpha(X, A)\cap C_0(X, A).$$

In the case where, $A=\Bbb{C}$, we write $\hbox{Lip}_\alpha(X, \Bbb{C})=\hbox{Lip}_\alpha X$, $\hbox{lip}_\alpha(X, \Bbb{C})=\hbox{lip}_\alpha X$ and $\hbox{lip}^\circ_\alpha (X, \Bbb{C})=\hbox{lip}^\circ_\alpha X$. For an extensive study of Lipschitz spaces see \cite{bcd, d, s2, s3, s4}.

Let $G$ be a locally compact group with identity element $e$ and weight function $\omega$, i.e., a continuous function $\omega: G\rightarrow [1, \infty)$ such that $\omega(e)=1$ and $\omega(xy)\leq\omega(x)\omega(y)$ for all $x, y\in G$. Let $L^1(G, \omega)$ and $M(G, \omega)$ be the Beurling and weighted measure algebra of $G$, respectively; for more details see \cite{dl, r0}.  We put $$L^1(G, 1)=L^1(G)\quad\hbox{and}\quad M(G, 1)=M(G).$$
For a Banach algebra $A$, let $L^1(G, A)$ be the Banach algebra of $A-$valued integrable functions on $G$. It is well-known that $L^1(G, A)$ isometrically isomorphic to the projective tensor product $L^1(G)\hat{\otimes} A$; see for example Proposition 1.5.4 in \cite{k}.

This paper is organized as follow. In Section 2, we give some results on homological properties of Banach algebras. For example, we prove that biprojectivity, contractibility and biflatness preserve under retractions. Also, if $A^{**}$ is contractible, then $A$ is contractible. Then, we show that the character space of a commutative contractible Banach algebra is finite. Section 3 is devoted to study homological properties of vector-valued Banach algebras. To this end, let $A$ be a Banach algebra such that $\Delta(A)$, the character space of $A$, is non-empty. Let also  $\frak{F}(X, A)$ be one of the Banach algebras $C_b(X, A)$, $C_0(X, A)$, $\hbox{Lip}_\alpha(X, A)$ or $\hbox{lip}_\alpha(X, A)$. We show that biprojectivity of $\frak{F}(X, A)$ implies discreteness of $X$ and $\Delta(A)$. We also prove that $\frak{F}(X, A)$ is contractible if and only if $X$ is finite and $A$ is contractible. Then, for commutative Banach algebra $A$, we prove that $\frak{F}(X, A)$ is contractible if and only if $A$ is a $C^*-$algebra and both $X$ and $\Delta(A)$ are finite. A special case of this result shows that $\hbox{lip}_\alpha^0X$ is contractible if and only if $X$ is finite. In Section 4, we investigate homological properties weighted Banach algebras and show that if $\frak{X}$ is a Banach algebra containing $C_0(G, 1/\omega)$, then contractibility and biprojectivity of $\frak{X}^*$ are equivalent to finitness of $G$, however, biflatness of $\frak{X}$ implies discreteness of $G$ and boundedness of the weight function $\omega^*$ defined by $\omega^*(x)=\omega(x)\omega(x^{-1})$ for all $x\in G$.
We also prove that biprojectivity, contractibility and biflatness of the Banach algebras $M(G, \omega)^{**}$ and $ L^1(G, \omega)^{**}$ are equivalent and occur only if $G$ is finite. Finally, we prove that  $L^1(G, A)$ is contractible if and only if $M(G)\hat{\otimes}A$ is contractible; or equivalently, $G$ is finite and $A$ is contractible.

\section{\normalsize\bf Homological and cohomological properties of Banach algebras}

Curtis and Loy \cite{cl} and Helemskii \cite{hel, hel1, h3, h4, hel2} investigated the relation between the homological properties and amenability and obtained the following result.

\begin{theorem}\label{homology} Let $A$ be a non-zero Banach algebra. Then the following statements hold.

\emph{(i)} Every $C^*-$ algebra is biflat.

\emph{(ii)} If $A$ is biflat, then $A$ is weakly amenable.

\emph{(iii)} If $A$ is biprojective, then $A$ is biflat.

\emph{(iv)} $A$ is contractible if and only if $A$ is unital and biprojective.

\emph{(v)} If $A$ is contractible, then $A$ is amenable.

\emph{(vi)} Commutative Banach algebra $A$ is contractible if and only if $A$ is finite dimensional.

\emph{(vii)} $A$ is amenable if and only if $A$ has a bounded approximate identity and $A$ is biflat.

\emph{(viii)} If $A$ is a commutative biprojective Banach algebra, then $\Delta(A)$ is discrete. Furthermore, if $A$ is a commutative $C^*-$algebra, then the converse remain valid.
\end{theorem}

\begin{example}{\rm Let $E$ be a Banach space with $\hbox{dim}(E)>1$ and let $\varphi\in E^*$ be non-zero. Then $E$ with the following multiplication is a Banach algebra.
$$
a\cdot b:=\varphi(b)a.
$$
It is clear that $A:=(E, \cdot)$ is biprojective and so $A$ is biflat. Also, $\Delta (A)=\{\varphi\}$. Note that $A$ does not have bounded approximate identity. Hence $A$ is neither contractible nor amenable. Thus the assumption that $A$ has a bounded approximate identity in Theorem \ref{homology} (iv) and (vii) can not drop.
}
\end{example}

Let $A$ and $B$ be Banach algebras. Then continuous homomorphism $\phi: A\rightarrow B$ is called a \emph{retraction} if there exists a continuous homomorphism $\psi: B\rightarrow A$ such that $\phi\circ\psi=id_{B}$, the identity map on $B$. 
Also, a linear map $L$ from $A$ into Banach left $A-$module $E$ is called a \emph{left multiplier} if $L(ab)=L(a)b$ for all $a, b\in A$. Similarly, one can define the notion right multiplier.  Pirkovskii and Selivanov \cite{psel} proved that $A$ is biprojective (respectively, biflat) if and only if every continuous derivation $D$ from $A$ into Banach $A-$bimodule $E$ (respectively, $E^*$), there exist left multiplier $L$ and right multiplier $R$ from $A$ into $E$ (respectively, $E^*$) such that $D=R-L$. 

\begin{theorem}\label{kilu} Let $A$ and $B$ be Banach algebras and $\phi: A\rightarrow B$ be a retraction. Then the following statements hold.

\emph{(i)} If $A$ is biprojective, then $B$ is biprojective.

\emph{(ii)} If $A$ is contractible, then $B$ is contractible.

\emph{(iii)} If $A$ is biflat, then $B$ is biflat.
\end{theorem}
{\it Proof.} First note that if $E$ is a Banach $B-$bimodule, then $E$ is a Banach $A-$bimodule with the following actions.
$$
a\cdot x=\phi(a)x\quad\hbox{and}\quad x\cdot a=x\phi(a)
$$
for all $a\in A$ and $x\in E$. Now, let $D: B\rightarrow E$ be a continuous derivation. Then the linear operator $\bar{D}: A\rightarrow E$ defined by $\bar{D}=D\circ\phi$ is a continuous derivation; indeed, for every $a_1, a_2\in A$, we have
\begin{eqnarray*}
\bar{D}(a_1a_2)&=&D(\phi(a_1)\phi(a_2))\\
&=&D(\phi(a_1))\phi(a_2)+\phi(a_1)D(\phi(a_2))\\
&=&\bar{D}(a_1)\cdot a_2+a_1\cdot\bar{D}(a_2).
\end{eqnarray*}
Let $A$ be biprojective. Then there exist left multiplier $L$ and right multiplier $R$ from $A$ into $E$ such that $\bar{D}=R-L$. Since $\phi$ is a retraction, there exists a bounded homomorphism $\psi: B\rightarrow A$ with $\phi\circ\psi=id_B$. Define left multiplier and right multiplier $\tilde{L}$ and $\tilde{R}$ by
$$
\tilde{L}=L\circ\psi\quad\hbox{and}\quad\tilde{R}=R\circ\psi.
$$
Then for every $b\in B$ we have
\begin{eqnarray*}
D(b)&=&D(\phi\psi(b))=\bar{D}(\psi(b))\\
&=&R\psi(b)-L\psi(b)\\
&=&\tilde{L}(b)-\tilde{R}(b).
\end{eqnarray*}
Therefore, $B$ is biprojective. That is, (i) holds. For (ii), note that if $A$ is unital, then $B$ is unital. This together with (i) proves (ii). The statement (iii) proves similar to (i).
$\hfill\square$

\begin{corollary} Let $A$ and $B$ be Banach algebras and $\phi: A\rightarrow B$ be an isomorphism and a homeomorphism. Then $\phi$ preserves biprojectivity, contractibility and biflatness.
\end{corollary}

 It is natural to ask whether continuous homomorphisms with dense range preserve biprojectivity, contractibility and biflatness. Mewomo \cite{mo} proved that the answer is affirmative for contractibility. But, the following example shows that the question has a negative answer for biprojectivity and biflatness. 

\begin{example}\label{nu}{\rm Let $E$ be a Banach space. Then $E\hat{\otimes} E^*$ with the following multiplication is a Banach algebra.
$$
(x_1\otimes f_1)\cdot(x_2\otimes f_2)=f_1(x_2)x_1\otimes f_2\quad\quad(x_1, x_2\in E, f_1, f_2\in E^*)
$$
Set $A:=(E\hat{\otimes} E^*, \cdot)$. Define the mapping $Q: A\hat{\otimes}A\rightarrow A$ by $$Q(x\otimes f)\otimes(y\otimes g)=f( y)x\otimes g.$$ Choose $y\in G$ such that $f(y)\neq 0$. Define the function $S: A\rightarrow A\hat{\otimes}A$ by
$$
S(x\otimes f)=( x\otimes f)\otimes( y\otimes 1/f(y)f).
$$
Then for every $x\in E$ and $f\in E^*$, we have
\begin{eqnarray*}
Q\circ S(x\otimes f)&=&
Q(x\otimes f)\otimes( y\otimes 1/f(y)f)\\
&=&
f(y)x\otimes 1/f(y)f\\
&=&x\otimes f .
\end{eqnarray*}
This shows that $Q\circ S=id_A$. It follows that $A$ is biprojective and hence it is biflat.
We define the function $\phi: E\hat{\otimes} E^*\rightarrow B(E)$ by $$\phi(a\otimes f)(x)= f(x) a.$$ If $X_i=x_i\otimes f_i$, for $i=1, 2$, then for every $x\in E$, we have
\begin{eqnarray*}
\phi(X_1)\circ \phi(X_2)(x)&=&\phi(X_1)[(f_2(x) x_2)]\\
&=&f_1( f_2(x)x_2)x_1\\
&=&f_2(x)f_1(x_2)x_1=\phi(X_1\cdot X_2)(x)
\end{eqnarray*}
Hence $\phi$ is a continuous homomorphism. Usually, the range of $\phi$ is called the \emph{space of }\emph{nuclear operators} on $E$ and is denoted by ${\cal N}(E)$. So $\phi$ is a continuous epimorphism from $A$ onto ${\cal N}(E)$.
It is well-known from \cite{gro3} that ${\cal N}(E)$ is weakly amenable if and only if $$\hbox{dim}(\hbox{ker} \phi)\leq 1.$$ So if $\hbox{dim}(\hbox{ker} \phi)\geq 2$, then ${\cal N}(E)$ is not weakly amenable and thus it is neither biprojective nor biflat. Therefore, continuous epimorphisms don't preserve biprojectivity and biflatness.

On can prove that $(\hbox{ker}\phi)^2=\{0\}$ and so it is not essential. This shows that $\hbox{ker}\phi$ is not weakly amenable. Note that $\hbox{ker}\phi$ is an ideal of Banach algebra $E\hat{\otimes} E^*$. So closed ideals of a biprojective (respectively, biflat) algebras is not  biprojective (respectively, biflat), in general.
}
\end{example}

\begin{proposition}\label{biprojective} Let $A_1$ and $A_2$ be Banach algebras. Then the following statements hold.

\emph{(i)} $A_1\oplus A_2$ is biprojective if and only if $A_1$ and $A_2$ are biprojective.

\emph{(ii)}  $A_1\oplus A_2$ is contractible if and only if $A_1$ and $A_2$ are contractible.

\emph{(iii)} $A_1\oplus A_2$ is biflat if and only if $A_1$ and $A_2$ are biflat.
\end{proposition}
{\it Proof.} (i) For $i=1, 2$, let $\pi_i: A_1\oplus A_2\rightarrow A_i$ and $\iota_i: A_i\rightarrow A_1\oplus A_2$ be canonical projection and injection, respectively. Then $\pi_i\circ\iota_i= id_{A_i}$, the identity map on $A_i$. Hence $\pi_i$ is a retraction of $A_1\oplus A_2$. So, if $A_1\oplus A_2$ is biprojective, then $A_i$ is biprojective. The converse is clear.

(ii) This follows from (i) and Theorem \ref{homology} (iv).

(iii) This is proved by the argument used in (i).$\hfill\square$\\

Ramsden \cite{ram1} showed that if $A_1$ and $A_2$ are biprojective (biflat), then $A_1\hat{\otimes}A_2$ is biprojective (biflat). He also proved that if $A_1$ is unital and $A_2$ containing a non-zero idempotent element, then biprojectivity (biflatness) of $A_1\hat{\otimes}A_2$ implies biprojectivity (biflatness) of $A_1$. 

\begin{theorem}\label{moka} Let $A_1$ and $A_2$ be Banach algebras. Then the following statements hold.

\emph{(i)} If $A_1$ and $A_2$ are contractible, then $A_1\hat{\otimes}A_2$ is contractible. The converse is true if $\Delta(A_i)$ is a nonempty set, for $i=1, 2$.

\emph{(ii)} If $A_2$ is unital and $\Delta(A_2)$ is a nonempty set, then biprojectivity \emph{(}biflatness\emph{)} of $A_1\hat{\otimes}A_2$ implies biprojectivity \emph{(}biflatness\emph{)}  of $A_1$.

\emph{(iii)} If $A_1$ is unital and $\Delta(A_1)$ is a nonempty set, then biprojectivity \emph{(}biflatness\emph{)}  of $A_1\hat{\otimes}A_2$ implies biprojectivity \emph{(}biflatness\emph{)}  of $A_2$.
\end{theorem}
{\it Proof.} Let $A_1$ and $A_2$ be contractible. Then $A_1$ and $A_2$ are unital and biprojective. Thus $A_1\hat{\otimes}A_2$ is unital and biprojective. So $A_1\hat{\otimes}A_2$ is contractible. 

Let $\varphi_2\in\Delta(A_2)$. Define the mapping $F: A_1\hat{\otimes} A_2\rightarrow A_1$ by $F(a_1\otimes a_2)=\varphi_{2}(a_{2})a_1$. Then for every $a_1, x_1\in A_1$ and $a_2, x_2\in A_2$ we have
\begin{eqnarray*}
F((a_1\otimes a_2)(x_1\otimes x_2))&=&
F(a_1x_1\otimes a_2x_2)\\
&=&
\varphi_{2}(a_{2}x_{2})a_1x_1\\
&=&
\varphi_{2}(a_{2})\varphi_{2}(x_{2}) a_1x_1\\
&=&
F(a_1\otimes a_2)F( x_1\otimes x_2).
\end{eqnarray*}
Choose $a_{2}\in A_{2}$ with $\varphi_{2}(a_{2})=1$. If $a_1\in A_1$, then $F(a_1\otimes a_2)=a_1$. Thus $F$ is a continuous epimorphism.
Hence contractibility of $A_1$ follows from the fact that every continuous epimorphism of Banach algebras is contractibility preserving.
Similarly, $A_2$ is contractible. So (i) holds. 

Assume that $A_2$ is unital and $\Delta(A_2)$ is a nonempty set. Define the function $\psi: A_1\rightarrow A_1\hat{\otimes}A_2$ by $\psi(a_1)=a_1\otimes 1_{A_2}$. Then $\psi$ is a homomorphism and $F\circ\psi=id_{A_1}$. Thus $F$ is a retraction. Therefore, (ii) holds. Similarly, (iii) is proved.$\hfill\square$\\

For a Banach algebra $A$, Moslehian and Niknam \cite{mnik} showed that biflatness of $A^{**}$ implies biflatness of $A$. Essmaili et al. \cite{esr} established that if $A$ is an ideal of $A^{**}$, then biprojectivity of $A^{**}$ forces $A$ to be biprojective. In the next result, we give relation between contractibility of $A^{**}$ and $A$. To this end, let us recall that a derivation $D$ from Banach algebra $A$ into Banach $A-$bimodule $E$ is called \emph{weakly approximately inner} if there exists a bounded net $(x_\alpha)$ in $E$ such that for every $a\in A$
$$
D(a)=w-\lim_\alpha\;x_\alpha.a-a.x_\alpha.
$$

\begin{theorem}\label{ek} Let $A$ be a Banach algebra. If $A^{**}$ is contractible, then $A$ is contractible.
\end{theorem}
{\it Proof.} Let $D$ be a continuous derivation from $A$ into Banach $A-$bimodule $E$. Then $\frak{B}=A\oplus E$ is a Banach algebra with the product 
$$
(a, x)(b, y)=(ab, ay+xb)
$$
and the norm $\|(a, x)\|=\|a\|+\|x\|$. Also, the mapping $\theta: A\rightarrow A\oplus E$ defined by $\theta(a)=(a, D(a))$ is a continuous homomorphism and $\theta^{**}(n)=(n, D^{**}(n))$ for all $n\in A^{**}$. Theses facts show that $D^{**}$ is a continuous derivation; see Lemma 2.2 in \cite{gstud}. Since $A^{**}$ is contractible, $A^{**}$ is amenable and so $D^{**}$ is inner. Hence there exists $\Phi\in E^{**}$ such that for every $m\in A^{**}$ and $f\in E^*$, we have
$$
\langle D^{**}(m), f\rangle=\langle \Phi\cdot m-m\cdot \Phi, f\rangle=\langle \Phi, m\cdot f-f\cdot m\rangle.
$$
Assume that $(x_\alpha)$ be a bounded net in $E$ such that $x_\alpha\rightarrow \Phi$ in the weak$^*$ topology of $E^{**}$. For every $a\in A$, we have
\begin{eqnarray*}
\langle D(a), f\rangle&=&
\langle D^{**}(a), f\rangle\\
&=&\Phi( a.f)-\Phi(f.a)\\
&=&
\lim_\alpha(a.f(x_\alpha)-f.a(x_\alpha))\\
&=&
\lim_\alpha\langle f, x_\alpha.a-a.x_\alpha\rangle.
\end{eqnarray*}
This implies that
$$
D(a)=w-\lim_\alpha\;x_\alpha.a-a.x_\alpha.
$$
Hence $D$ is weakly approximately inner and therefore it is inner. That is, $A$ is contractible.$\hfill\square$\\

 Essmaili et al. \cite{esr} raised the question of whether there exists a Banach algebra $A$ such that $A^{**}$ is biprojective, but $A$ is not biprojective. The next example gives a positive answer to the question. This example also shows that the converse of Theorem \ref{ek} is not true, in general. 

\begin{example}{\rm Let $A=c_0$. Then $A$ is biprojective. If $A^{**}=\ell^\infty$ is biprojective, then it is commutative and contractible. So $A^{**}$ is finite dimensional, a contradiction. Thus $A^{**}$ is not biprojective. 
Note that $A^\sharp$, the unitization of $A$, is contractible, however, $(A^\sharp)^{**}$ is a commutative infinite dimensional Banach algebra. So $(A^\sharp)^{**}$ is not contractible.}
\end{example}

It is well-known that the character space of unital, commutative Banach algebras is weak$^*$ compact; see for example Theorem 2.2.3 in \cite{k}.
Here, we give a simple proof for it. Note that an application of Shilov's idempotent theorem shows that if the character space of a semisimple commutative Banach algebra is compact, then the Banach algebra is unital; see Theorems 3.5.1 and 3.5.5 in \cite{k}.

\begin{lemma}\label{del} Let $A$ be a unital, commutative Banach algebra. Then $\Delta(A)$ is weak$^*$ compact.
\end{lemma}
{\it Proof.} Let $(\varphi_\alpha)$ be a net in $\Delta(A)$ such that $\varphi_\alpha\rightarrow\varphi\in A^*$ in the weak$^*$ topology of $A^*$.  Then
$$
\varphi(1_A)=\lim_\alpha\varphi_\alpha(1_A)=1.
$$
So $\varphi$ is non-zero. If $f, g\in A$, then
\begin{eqnarray*}
\varphi(fg)&=&\lim_\alpha\varphi_\alpha(fg)=\varphi_\alpha(fg)\\
&=&\varphi_\alpha(f)\varphi_\alpha(g)=\varphi(f)\varphi(g).
\end{eqnarray*}
Hence $\varphi$ is multiplicative. Since $A$ is commutative, by \cite{k}, $\varphi$ is continuous. It follows that $\varphi\in\Delta(A)$. Thus $\Delta(A)$ is a weak$^*$ closed subspace of the compact unit ball of $A^*$. Therefore, $\Delta(A)$ is weak$^*$ compact.$\hfill\square$\\

We finish this section with the following result.

\begin{theorem} Let $A$ be a  commutative Banach algebra. Then the following assertions are equivalent

\emph{(a)} $A$ is contractible.

\emph{(b)}$A$ is finite dimensional.

\emph{(c)} $A$ is a $C^*-$algebra and $\Delta(A)$ is finite.
\end{theorem}
{\it Proof.} By Theorem \ref{homology} the statements (a) and (b) are equivalent. Let $A$ be contractible. Then $A$ is unital and biprojective. By  Theorem \ref{homology} (viii) and Lemma \ref{del}, $\Delta(A)$ is finite. Since $A$ is commutative, $A$ is a $C^*-$algebra; see Corollary 4.1.3 in \cite{run}. Thus (a)$\Rightarrow$(c). Finally, the implication (c)$\Rightarrow$(a) follows from Theorem \ref{homology} (iv) and (vii).$\hfill\square$

\section{\normalsize\bf Homological and cohomological properties of vector-valued function algebras}

Let $X$ be a locally compact space and $A$ be a Banach algebra.  For $x\in X$ and $\varphi\in\Delta(A)$, we define the multiplicative linear functional $x\otimes\varphi: X^A\rightarrow{\Bbb C}$ by
$$
x\otimes\varphi(f):=\varphi(f(x))
$$
for all $f\in X^A$.
 We set
$$
X\otimes\Delta(A)=\{x\otimes\varphi: x\in X, \varphi\in\Delta(A)\}.
$$

\begin{proposition}\label{book1} Let $(X, d)$ be a locally compact Hausdorff space and $A$ be a Banach algebra with $\Delta(A)\neq\emptyset$. Then the function $\tilde{\Gamma}: (x, \varphi)\mapsto x\otimes\varphi$ from $X\times\Delta(A)$ onto $X\otimes\Delta(A)$ is separately continuous and bijective.
\end{proposition}
{\it Proof.} Let $x_0$ and $x_1$ be distinct points of $X$. Then there exists a function $f\in C_c(X)$ such that $f(x_0)=1$ and $f(x_1)=0$. Now, let $\varphi$ and $\psi$ be distinct elements of $\Delta(A)$. Choose $a\in A$ such that $\varphi(a)$ is non-zero. Define the function $g: X\rightarrow A$ by $g(x)=af(x)$. Then $g(x_0)=a$ and $g(x_1)=0$. So $\varphi(g(x_0))\neq\psi(g(x_1))$ . Hence, $x_0\otimes\varphi\neq x_1\otimes\psi$. Thus $\tilde{\Gamma}$ is injective. It is clear that $\tilde{\Gamma}$ is surjective and separately continuous.$\hfill\square$\\

Let $\frak{A}(X)$ be one of the Banach algebras $C_b(X)$, $\hbox{Lip}_\alpha X$ or $\hbox{lip}_\alpha X$, where $\alpha>0$. Let also $\frak{F}(X)$ be $\frak{A}(X)$ or $C_0(X)$. Note that in the case where $\frak{F}(X)$ is $C_b(X)$ or $C_0(X)$, the set $X$ is a locally compact space and the other cases, $X$ is a metric space. For a Banach algebra $A$, let $\frak{A}(X, A)$ and $\frak{F}(X, A)$ be the vector-valued Banach algebras corresponding to $\frak{A}(X)$ and $\frak{F}(X)$, respectively.

\begin{lemma}\label{re} Let $A$ be a Banach algebra. Then the following statements hold.

\emph{(i)} If $\varphi\in\Delta(A)$, then $\frak{F}(X)=\varphi\otimes\frak{F}(X, A)$, where $\varphi\otimes f(x)=\varphi(f(x))$ for all $f\in \frak{F}(X)$ and $x\in X$.

\emph{(ii)} If $x\in X$, then $x\otimes\frak{A}(X, A)=A$, where $x\otimes g=g(x)$ for all $g\in \frak{A}(X, A)$.

\emph{(iii)}  For every $\varphi\in\Delta(A)$, the mapping $f\mapsto \varphi\otimes f$ is a continuous epimorphism from $\frak{F}(X, A)$ onto $\frak{F}(X)$.

\emph{(iv)} For every $x\in X$, the mapping $f\mapsto x\otimes f$ is a continuous epimorphism from $\frak{A}(X, A)$ onto $A$.

\emph{(v)} $\frak{A}(X, A)$ is unital if and only if $A$ is unital.

\emph{(vi)} $\frak{A}(X, A)$ is commutative if and only if $A$ is commutative.

\end{lemma}
{\it Proof.} Let $\varphi\in\Delta(A)$ and $g\in \frak{F}(X)$. Choose $a\in A$ such that $\varphi(a)=1$. Then $ag\in \frak{F}(X, A)$ and
$$
\varphi\otimes(ag)(x)=\varphi(ag(x))=g(x)\varphi(a)=g(x)
$$
for all $x\in X$. So $g=\varphi\otimes(ag)$. Thus $\frak{F}(X)$ contains in $\varphi\otimes\frak{F}(X, A)$. It is easy to see that $\varphi\otimes f\in \frak{F}(X)$ for all $f\in \frak{F}(X, A)$. So (i) holds.

Let $a\in A$. Define the vector-valued function $f\in \frak{A}(X, A)$ by $f(x)=a$ for all $x\in X$. Then $a=x\otimes f$. This shows that (ii) holds.

The statements (iii) and (iv) follow from (i) and (ii), respectively. For (v), let $u$ be an identity element of $\frak{A}(X, A)$. Fix $x_0\in X$. Then for every $a\in A$, we have
$$
au(x_0)=x_0\otimes (au)=x_0\otimes a=a.
$$
Hence $u(x_0)$ is an identity of $A$. Thus (v) holds. The proof of statement (vi) is rutin.$\hfill\square$

\begin{proposition} Let $A$ be a commutative Banach algebra with $\Delta(A)\neq\emptyset$ and $\alpha>0$.  If $\frak{F}(X, A)$ is biprojective, then $X$ and $\Delta(A)$ are discrete.
%
%
%
%
\end{proposition}
{\it Proof.}  Let $\frak{F}(X, A)$ be biprojective. In view of Theorem \ref{del}, $\Delta(\frak{F}(X, A))$ is discrete. Since $X\otimes\Delta(A)$ is a subset of $\Delta(\frak{F}(X, A))$,  it follows that $X$ and $\Delta(A)$ are discrete.$\hfill\square$\\

It is well-known that if $X$ is a locally compact space, then $\Delta(C_0(X, A))$ and $X\otimes\Delta(A)$ are homeomorphic.

\begin{theorem}\label{zmah1} Let $A$ be a commutative Banach algebra with $\Delta(A)\neq\emptyset$. Then the following assertions are equivalent.

\emph{(a)} $C_b(X, A)^{**}$ is contractible.

\emph{(b)} $C_b(X, A)$ is contractible.

\emph{(c)} $C_0(X, A)^{**}$ is contractible.

\emph{(d)} $C_0(X, A)$ is contractible.

\emph{(e)} $X$ is finite and $A$ is finite dimensional.

\emph{(f)} $A$ is a $C^*-$algebra and both $X$ and $\Delta(A)$ are finite.
\end{theorem}
{\it Proof.} 
It follows from Theorem \ref{ek} that (a)$\Rightarrow$(b) and (c)$\Rightarrow$(d). Let $\frak{F}(X, A)$ be contractible. Then $\Delta(\frak{F}(X, A))$ is finite. Since $X\otimes\Delta(A)$ is a subset of $\Delta(\frak{F}(X, A))$, it follows that $X\otimes\Delta(A)$ is finite. So $X$ is finite.
In view of Lemma \ref{re} and Proposition 2.4.2 in \cite{mo}, $A$ is unital, commutative and contractible. Hence $A$ is finite dimensional. Thus (b)$\Rightarrow$ (e) and (d)$\Rightarrow$(e).

If $A$ is finite dimensional, then there exists $n\in{\Bbb N}$ such that $A$ and ${\Bbb C}^n$ are isomorphism. Thus $A$ is a contractible $C^*-$algebra and so $\Delta(A)$ is finite. Hence (e)$\Rightarrow$(f).

Let $X$ be finite. Then $\frak{F}(X, A)=C_0(X, A)$ and so
$$
\Delta(\frak{F}(X, A))=X\otimes\Delta(A).
$$
If $A$ is a $C^*-$algebra and $\Delta(A)$ is finite, then $\frak{F}(X, A)$ is a commutative $C^*-$algebra and its character space is finite. By Theorem \ref{del}, $\frak{F}(X, A)$ is contractible. That is, (f)$\Rightarrow$(b) and (f)$\Rightarrow$(d).

Finally, assume that $X$ is finite and $A$ is finite dimensional. Then there exists $n\in {\Bbb N}$ such that $\frak{B}(X, A)^{**}$ is a finite number of copies ${\Bbb C}\hat{\otimes}{\Bbb C}^n$. Hence $\frak{F}(X, A)^{**}$ is contractible. That is, (e)$\Rightarrow$(a) and (e)$\Rightarrow$(c).
$\hfill\square$\\

The argument used in the proof of Theorem \ref{zmah1} shows the next result holds.

\begin{theorem}\label{rozi} Let $A$ be a commutative Banach algebra with $\Delta(A)\neq\emptyset$. Then the following assertions are equivalent.

\emph{(a)} $\emph{Lip}_\alpha(X, A)^{**}$ is contractible.

\emph{(b)} $\emph{Lip}_\alpha(X, A)$ is contractible.

\emph{(c)} $\emph{lip}_\alpha(X, A)^{**}$ is contractible.

\emph{(d)} $\emph{lip}_\alpha(X, A)$ is contractible.

\emph{(e)} $X$ is finite and $A$ is finite dimensional.

\emph{(f)} $A$ is a $C^*-$algebra and both $X$ and $\Delta(A)$ are finite.
\end{theorem}

\begin{corollary}\label{zmahh} Let $(X, d)$ be a locally compact metric space, $A$ be a Banach algebra and $\alpha>0$. Then the following statements hold.

\emph{(i)}  $\emph{\hbox{lip}}^0_\alpha(X, A)$ is contractible if and only if $X$ is finite and $A$ is finite dimensional.

\emph{(ii)} $\emph{\hbox{lip}}^0_\alpha X$ is contractible if and only if $X$ is finite.
\end{corollary}
{\it Proof.} Let $\hbox{lip}^0_\alpha(X, A)$ be contractible. Hence it is unital and so $C_0(X, A)$ is unital. Thus $X$ is compact. Now, apply Theorem \ref{rozi} to conclude that (i) holds. The statements (ii) follows at once from (i).$\hfill\square$


\begin{theorem}\label{zmah} Let $A$ be a Banach algebra with $\Delta(A)\neq\emptyset$. Then the following assertions are equivalent.

\emph{(a)} $C_b(X, A)$ \emph{(}respectively, $C_b(X, A)^{**}$\emph{)} is contractible.

\emph{(b)} $C_0(X, A)$ \emph{(}respectively, $C_0(X, A)^{**}$\emph{)} is contractible.

\emph{(c)} $X$ is finite and $A$ \emph{(}respectively, $A^{**}$\emph{)} is contractible.
\end{theorem}
{\it Proof.} Let $\frak{F}(X, A)$ be contractible. Then  $\frak{A}(X, A)$ or $C_0(X, A)$ is contractible.  Assume that $\frak{A}(X, A)$ is contractible. Then Proposition 2.4.2 in \cite{mo} and Lemma \ref{re} show that $\frak{A}(X)$ and $A$ are contractible. In view of Theorem \ref{zmah1}, $X$ is finite. Now, if $C_0(X, A)$ is contractible, then it is unital. So $X$ is compact and hence $C_0(X, A)=C_b(X, A)$. This implies that $X$ is finite and $A$ is contractible.

For the converse, note that if $X$ is finite, then $\frak{F}(X, A)$ is a finite number of copies ${\Bbb C}\hat{\otimes}A$. So, if $A$ is contractible, then by Theorem \ref{moka},  $\frak{A}(X, A)$ be contractible.$\hfill\square$\\

The proof of the following result is similar to Theorem \ref{zmah}, and so we omit it.

\begin{theorem} Let $A$ be a Banach algebra with $\Delta(A)\neq\emptyset$. Then the following assertions are equivalent.

\emph{(a)} $\emph{Lip}_\alpha(X, A)$ \emph{(}respectively, $\emph{Lip}_\alpha(X, A)^{**}$\emph{)} is contractible.

\emph{(b)} $\emph{lip}_\alpha(X, A)$ \emph{(}respectively, $\emph{lip}_\alpha(X, A)^{**}$\emph{)} is contractible.

\emph{(c)} $X$ is finite and $A$ \emph{(}respectively, $A^{**}$\emph{)} is contractible.
\end{theorem}

Theorem \ref{zmah} together with Theorem \ref{homology} and Lemma \ref{re}  proves the following results.

\begin{corollary} Let $A$ be a unital Banach algebra with $\Delta(A)\neq\emptyset$. Then the following assertions are equivalent.

\emph{(a)} $C_b(X, A)$ \emph{(}respectively, $C_b(X, A)^{**}$\emph{)} is biprojective.

\emph{(b)} $C_b(X, A)$ \emph{(}respectively, $C_b(X, A)^{**}$\emph{)} is contractible.

\emph{(c)} $X$ is finite and $A$ \emph{(}respectively, $A^{**}$\emph{)} is contractible.
\end{corollary}

\begin{corollary} Let $A$ be a unital Banach algebra with $\Delta(A)\neq\emptyset$. Then the following assertions are equivalent.

\emph{(a)} $\emph{Lip}_\alpha(X, A)$ \emph{(}respectively, $\emph{Lip}_\alpha(X, A)^{**}$\emph{)} is biprojective.

\emph{(b)} $\emph{Lip}_\alpha(X, A)$ \emph{(}respectively, $\emph{Lip}_\alpha(X, A)^{**}$\emph{)} is contractible.

\emph{(c)} $\emph{lip}_\alpha(X, A)$ \emph{(}respectively, $\emph{lip}_\alpha(X, A)^{**}$\emph{)} is biprojective.

\emph{(d)} $\emph{lip}_\alpha(X, A)$ \emph{(}respectively, $\emph{lip}_\alpha(X, A)^{**}$\emph{)} is contractible.

\emph{(e)} $X$ is finite and $A$ \emph{(}respectively, $A^{**}$\emph{)} is contractible.
\end{corollary}


\section{\normalsize\bf Homological and cohomological properties of weighted group algebras}

Theorem \ref{homology} shows that there is a closed connection between cohomology and homology properties. This result enables us to find facts about homological  properties of group algebras and their second dual.

\begin{proposition} Let $\omega$ be a weight function on a locally compact group $G$. Then the following assertions are equivalent.

\emph{(a)} $L^1(G, \omega)$ is biflat.

\emph{(b)} $L^1(G, \omega)$ is  amenable.

\emph{(c)}  $G$ is amenable and $\omega^*$ is bounded
\end{proposition}
{\it Proof.} This follows from Theorem \ref{homology}  and the fact that $L^1(G, \omega)$ is amenable if and only if $G$ is amenable and $\omega^*$ is bounded; see \cite{gro3}.$\hfill\square$

Helmeski \cite{h3} showed that $L^1(G)$ is biprojective if and only if $G$ is compact. We prove this result for Beurling algebra $L^1(G, \omega)$.

\begin{proposition}\label{l1} Let $G$ be a locally compact group. Then the following statements hold.

\emph{(i)} $L^1(G, \omega)$ is biprojective if and only if $G$ is compact.

\emph{(ii)} $L^1(G, \omega)$ is contractible if and only if $G$ is finite.
\end{proposition}
{\it Proof.} (i) Let $L^1(G, \omega)$ be biprojective. Then  $L^1(G, \omega)$ is biflat and so it is amenable. In this case, it is known from \cite{w} that  the Banach algebras $L^1(G, \omega)$ and $L^1(G)$ are isomorphic and homeomorphic as Banach algebras. So Theorem \ref{kilu} shows that $L^1(G)$ is biprojective. Hence $G$ is compact. For the converse, we only remark that if $G$ is compact, then $L^1(G, \omega)=L^1(G)$.

(ii) Let $L^1(G, \omega)$ be contractible. Then $L^1(G, \omega)$ is unital and so $G$ is discrete. On the other hand, $L^1(G, \omega)$ is biprojective. By (i) $G$ is also compact. Therefore, $G$ is finite. Conversely, if $G$ is finite, then $L^1(G, \omega)$ is unital and by (i), it is biprojective. Thus $L^1(G, \omega)$ is contractible. $\hfill\square$\\

In the sequel, we investigate homological properties of weighted measure algebra $M(G, \omega)$.

\begin{proposition}\label{m} Let $\omega$ be a weight function on a locally compact group $G$. Then the following assertions are equivalent.

\emph{(a)} $M(G, \omega)$ is biflat.

\emph{(b)} $M(G, \omega)$ is  amenable.

\emph{(c)} $G$ is a discrete amenable group and $\omega^*$ is bounded.
\end{proposition}
{\it Proof.} Let $M(G, \omega)$ is biflat. In view of Proposition 2.4 , $L^1(G, \omega)$ is biflat and so $G$ is an amenable group and $\omega^*$ is bounded. On the other hand, $M(G, \omega)$ is weakly amenable. Hence $G$ is discrete; see Theorem 2.4 in \cite{mr5}. This together with  Theorem \ref{homology} proves the result.$\hfill\square$

\begin{proposition}\label{m} Let $\omega$ be a weight function on a locally compact group $G$. Then the following assertions are equivalent.

\emph{(a)} $M(G, \omega)$ is biprojective.

\emph{(b)} $M(G, \omega)$ is  contractible.

\emph{(c)} $G$ is finite.
\end{proposition}
{\it Proof.}  Since $M(G, \omega)$ is unital, (a) and (b) are equivalent. Let $M(G, \omega)$ be contractible. This implies that  $M(G, \omega)$ is weakly amenable. It follows from Theorem 2.4 in \cite{mr5} that $G$ is a discrete.
Thus $L^1(G, \omega)$ is contractible. By Proposition \ref{l1}, $G$ is finite.  That is (b) implies (c). Clearly, (c) implies (a).$\hfill\square$\\

Let  $C_0(G, 1/\omega)$ be the Banach algebra of all complex-valued  functions $f$ on $G$ such that  $f/\omega\in C_0(G)$.

\begin{theorem}\label{intro} Let $\omega$ be a weight function on a locally compact group $G$ and $\frak{X}$ be a Banach algebra containing $C_0(G, 1/\omega)$. Then the following assertions are equivalent.

\emph{(a)} $\frak{X}^*$ is biprojective.

\emph{(b)} $\frak{X}^*$ is contractible.

\emph{(c)} $G$ is finite.
\end{theorem}
{\it Proof.} Let $\frak{X}$ be a Banach algebra containing $C_0(G, 1/\omega)$. Then
$$
\frak{X}^*= M(G, \omega)\oplus C_0(G, 1/\omega)^\perp.
$$
If $\frak{X}^*$ is biprojective, then by Proposition \ref{biprojective}, $M(G, \omega)$ is biprojective. Hence $G$ is finite. So (a)$\Rightarrow$ (c). The implications (c)$\Rightarrow$(b)$\Rightarrow$(a) are clear.$\hfill\square$\\

\begin{theorem}\label{biflat} Let $\omega$ be a weight function on a locally compact group $G$ and $\frak{X}$ be a Banach algebra containing $C_0(G, 1/\omega)$. If $\frak{X}^*$ is biflat, then $G$ is a discrete amenable group and $\omega^*$ is bounded.
\end{theorem}
{\it Proof.} Let $\frak{X}^*$ be biflat. Since $M(G, \omega)$ is a direct summand of $\frak{X}^*$, it follows that $M(G, \omega)$ is biflat. Hence $G$ is a discrete amenable group and $\omega^*$ is bounded.$\hfill\square$\\

Let $C_b(G, 1/\omega)$ be the Banach algebra of all complex-valued  functions $f$ on $G$ such that $f/\omega\in C_b(G)$, and let  $ LUC(G, 1/\omega)$ be the Banach algebra of all $f\in C_b (G, 1/\omega)$ such that the map $x\mapsto\;_x(f/\omega)$ from $G$ into $C_b (G, 1/\omega)$ is norm continuous. Let also $\hbox{Ap}(G, 1/\omega)$ (respectively, $\hbox{Wap}(G, 1/\omega)$) be the set of all functions $f\in C_b (G, 1/\omega)$ such that  
$$
\{\frac{_xf}{\omega(x)\omega}: x\in G\}
$$
is relatively norm (respectively, weakly) compact in $C_b (G)$, where $_xf(y)=f(yx)$ for all $x, y\in G$.  Let $L_0^\infty(G, 1/\omega)^*$ and $M(G, \omega)_0^*$ be the Banach algebras defined as in \cite{mnr, mm}, respectively. It is easy to see that $C_0(G, 1/\omega)$ is a subspace of these Banach algebras. So, we may apply Theorems \ref{intro} and \ref{biflat} for them.  

\begin{theorem}\label{bi} Let $G$ be a locally compact group. Then the following assertions are equivalent.

\emph{(a)}  $L^1(G, \omega)^{**}$ \emph{(}respectively, $M(G, \omega)^{**}$\emph{)} is amenable.

\emph{(b)} $L^1(G, \omega)^{**}$ \emph{(}respectively, $M(G, \omega)^{**}$\emph{)}  is biflat.

\emph{(c)}  $L^1(G, \omega)^{**}$ \emph{(}respectively, $M(G, \omega)^{**}$\emph{)}  is biprojective.

\emph{(d)}  $L^1(G, \omega)^{**}$ \emph{(}respectively, $M(G, \omega)^{**}$\emph{)}  is contractible.

\emph{(e)}  $L^1(G, \omega)^{**}$ \emph{(}respectively, $M(G, \omega)^{**}$\emph{)} is weakly amenable.

\emph{(f)} $G$ is finite.
\end{theorem}
{\it Proof.}  In view of Theorem \ref{homology}, the statements (a)- (d) imply the statement (e). Also, by Theorem 3.1 in  \cite{mr5}, the statement (e) follows (f). So the proof will be completed if we note that the statement (f) implies the other statements.$\hfill\square$

\begin{remark}{\rm From Proposition \ref{l1} and Theorem\ref{bi} we see that if $G$ is an infinite compact group, then $L^1(G)$ is both biprojective and biflat, however, $L^1(G)^{**}$ is not neither biprojective nor bifalt. So biprojectivity (respectively, biflatness) of $A$ does not imply biprojectivity (respectively, biflatness) of $A^{**}$.}
\end{remark}

\begin{theorem}\label{bahar} Let $G$ be a locally compact group and $A$ be a Banach algebra with $\Delta(A)\neq\emptyset$. Then  $L^1(G, A)$ \emph{(}respectively, $L^1(G, A)^{**}$\emph{)} is contractible if and only if $G$ is finite and $A$ \emph{(}respectively, $A^{**}$ \emph{)} is contractible.
\end{theorem}
{\it Proof.} Let $\varphi\in\Delta(A)$. Then the mapping $f\mapsto\varphi\circ f$ is a continuous epimorphism from $L^1(G, A)$ onto $L^1(G)$. Also, if $x\in G$, then the mapping $f\mapsto f(x)$ is a continuous epimorphism from $L^1(G, A)$ onto $A$. Now, apply Proposition 2.4.2 in \cite{mo} and Proposition \ref{l1}.$\hfill\square$

\begin{theorem} Let $G$ be a locally compact group and $A$ be a unital, commutative Banach algebra. Then the following assertions are equivalent.

\emph{(a)} $L^1(G, A)$ is biflat.

\emph{(b)} $L^1(G, A)$ is amenable.

\emph{(c)} $G$ and $A$ are amenable.
\end{theorem}
{\it Proof.} Since $A$ is unital and $L^1(G)$ has a bounded approximate identity, (a) and (b) are equivalent. If (c) holds, then $L^1(G)$ and $A$ are amenable. So $L^1(G)\hat{\otimes}A$ is amenable. From this and the fact that $L^1(G, A)$ isometrically isomorphic to $L^1(G)\hat{\otimes}A$ we infer that (c) implies (b). By the proof of Theorem \ref{bahar}, the Banach algebras $L^1(G)$ and $A$ are continuous homomorphic images of $L^1(G, A)$. Hence amenability of $L^1(G, A)$ implies amenability of $G$ and $A$. That is, (b) implies (c).$\hfill\square$

\begin{theorem} Let $G$ be a locally compact group and $A$ be a Banach algebra with $\Delta(A)\neq\emptyset$.  If $A$ has a non-zero idempotent element, then the following statements hold.

\emph{(i)} $M(G)\hat{\otimes} A$ is biprojective if and only if $G$ is finite and $\ell^1(G)\hat{\otimes}A$ is biprojective.

\emph{(ii)} $M(G)\hat{\otimes} A$ is biflat if and only if $G$ is finite and $\ell^1(G)\hat{\otimes}A$ is biflat.
\end{theorem}
{\it Proof.} (i) Let $M(G)\hat{\otimes} A$ be biprojective. By \cite{ram1}, the Banach algebras $M(G)$ is biprojective. It follows from Theorem \ref{m} that $G$ is finite and so $M(G)\hat{\otimes} A= \ell^1(G)\hat{\otimes}A$ is biprojective. The converse is clear.

(ii) This proves similar to (i).$\hfill\square$

\begin{theorem} Let $G$ be a locally compact group and $A$ be a Banach algebra with $\Delta(A)\neq\emptyset$. Then $M(G)\hat{\otimes} A$ is contractible if and only if  $G$ is finite and $A$ is contractible.
\end{theorem}
{\it Proof.} Let $M(G)\hat{\otimes} A$ be contractible. By Theorem \ref{moka}, the Banach algebras $M(G)$ and $A$ are contractible. It follows from Theorem \ref{m} that $G$ is finite.$\hfill\square$

\vspace{2mm}

 {\footnotesize
\noindent {\bf Mohammad Javad Mehdipour}\\
Department of Mathematics,\\ Shiraz University of Technology,\\
Shiraz
71555-313, Iran\\ e-mail: mehdipour@.ac.ir\\
{\bf Ali Rejali}\\
Department of Pure Mathematics,\\ Faculty of Mathematics and Statistics,\\ University of Isfahan,\\
Isfahan
81746-73441, Iran\\ e-mail: rejali@sci.ui.ac.ir\\
\end{document}